\documentclass{amsart}

\usepackage{amssymb}

\newtheorem{thm}{Theorem}

\newtheorem{lem}[thm]{Lemma}
\newtheorem{cor}[thm]{Corollary}

\newtheorem{prop}[thm]{Proposition}

\theoremstyle{definition}
\newtheorem{defn}[thm]{Definition}

\newtheorem{say}[thm]{}

\newtheorem{rem}[thm]{Remark}          

\newtheorem*{ack}{Acknowledgments}      

\newtheorem{defn-thm}[thm]{Definition--Theorem}  
\newtheorem{defn-lem}[thm]{Definition--Lemma}  

\theoremstyle{remark}
\newtheorem{claim}[thm]{Claim}

\renewcommand{\o}[0]{{\mathcal O}} 
\newcommand{\z}[0]{{\mathbb Z}}
\newcommand{\n}[0]{{\mathbb N}}

\renewcommand{\a}[0]{{\mathbb A}}

\newcommand{\p}[0]{{\mathbb P}}

\newcommand{\q}[0]{{\mathbb Q}}

\newcommand{\qtq}[1]{\quad\mbox{#1}\quad}
\newcommand{\spec}[0]{\operatorname{Spec}}

\newcommand{\supp}[0]{\operatorname{Supp}}    
\newcommand{\red}[0]{\operatorname{red}}    
\newcommand{\codim}[0]{\operatorname{codim}}

\newcommand{\ext}[0]{\operatorname{Ext}}

\newcommand{\sing}[0]{\operatorname{Sing}}

\newcommand{\rdown}[1]{\lfloor{#1}\rfloor}

\newcommand{\onto}[0]{\twoheadrightarrow}

\newcommand{\simq}[0]{\sim_{\q}}

\newcommand{\depth}[0]{\operatorname{depth}} 
\newcommand{\tsum}[0]{\textstyle{\sum}}

\newcommand{\sext}[0]{\operatorname{\mathcal{E}\!\it{xt}}}    
 
\newcommand{\simql}[0]{\sim_{\q,loc}}

\newcommand{\nklt}[0]{\operatorname{nklt}}

\def\into{\DOTSB\lhook\joinrel\to}

\begin{document}
\bibliographystyle{amsalpha}

\title[Local Kawamata-Viehweg vanishing]{A local version of the \\ 
Kawamata-Viehweg vanishing theorem}
\author{J\'anos Koll\'ar}

\today

\maketitle

\hfill{\it To the memory of Eckart Viehweg}
\medskip

The, by now classical,  Kawamata--Viehweg vanishing theorem
\cite{kaw, vie}  says that {\em global} cohomologies vanish
for divisorial sheaves which are $\q$-linearly equivalent
to a divisor of the form $(\mbox{nef and big})+\Delta$.
In this note we prove that
  {\em local} cohomologies vanish
for divisorial sheaves which are $\q$-linearly equivalent
to a divisor of the form $\Delta$. 
If $X$ is a cone over a Fano variety, one can set up a perfect correspondence
between the global and local versions.

More generally,  we study the depth of various sheaves 
associated to a log canonical pair $(X,\Delta)$. 
The first significant result in this direction, due to \cite{elkik},
says that if $(X,0)$ is canonical then $X$ has rational
singularities. In particular, $\o_X$ is CM.
The proof has been simplified repeatedly 
in \cite{fujita85}, \cite[Sec.11]{k-pairs}
and \cite[5.22]{km-book}. Various generalizations 
for other divisorial sheaves and  to the log canonical case
were established in \cite[5.25]{km-book}, \cite{kovacs},  
\cite{ale-lim} and  \cite[Secs.4.2--3]{fujinobook}.

Here we prove a further generalization, which, I believe,
covers all the theorems about depth mentioned above.

We work with varieties over a field of characteristic 0.
For the basic definitions and for background material
see \cite{km-book}.

\begin{defn} Let $X$ be normal and $D_1, D_2$ two
$\q$-divisors. We say that $D_1$ is {\it locally $\q$-linearly equivalent}
to $D_2$, denoted by $D_1\simql D_2$,
 if $D_1-D_2$ is $\q$-Cartier.
The same definition works if $X$ is  not normal, as long as
none of the irreducible components of the $D_i$
is contained in $\sing X$. 

Note that this is indeed a local property.
That is, if $\{X_i: i\in I\}$ is an open cover of $X$ and
$D_1|_{X_i}\simql D_2|_{X_i}$ for every $i$, then $D_1\simql D_2$.
\end{defn}

The following can be viewed as a  local version of the 
Kawamata-Viehweg vanishing theorem.

\begin{thm}\label{d.num.eq.to.Delta.CM.thm}
 Let $(X,\Delta)$ be dlt,  $D$  a 
(not necessarily effective) $\z$-divisor 
and $\Delta'\leq \Delta $ an effective $\q$-divisor on
$X$ such that $D\simql \Delta'$.
Then $\o_X(-D)$ is CM.
\end{thm}

Here dlt is short for
divisorial log terminal \cite[2.37]{km-book} and
 CM  for Cohen--Macaulay. A sheaf $F$ is CM iff  
$H^i_x(X,F)=0$ for $i<\codim_Xx$ for every point $x\in X$
and 
$\depth_x F\geq j$ iff $H^i_x(X,F)=0$ for $i<j$; see
\cite[III.3.1]{sga2} or \cite[Exrcs.III.3.3--5]{hartsh}.

Examples illustrating the necessity of the assumptions are
given in (\ref{d.num.eq.to.Delta.CM.appl}.5--10).

The proof of  Theorem \ref{d.num.eq.to.Delta.CM.thm}
also works  in the complex analytic case.
(Normally one would expect that proofs of a local statement as above
automatically work for analytic spaces as well. However, many of the
papers cited above
 use global techniques, and some basic questions are still unsettled;
see, for instance, (\ref{d.num.eq.to.Delta.CM.appl}.3).)

\medskip

Weaker results hold for log canonical and semi log canonical pairs.
For basic definitions in the  semi log canonical case
(abbreviated as {\it slc})
see \cite[Sec.12]{k-etal} or \cite{fujinobook}.

\begin{thm}\label{main.slc.thm}
 Let $(X,\Delta)$ be semi log canonical and $x\in X$  
a point that is not a log canonical center of $(X,\Delta)$.
\begin{enumerate}
\item Let $D$ be a  $\z$-divisor such that
none of the irreducible components of  $D$
is contained in $\sing X$. Let
 $\Delta'\leq \Delta $ be an effective $\q$-divisor on
$X$ such that $D\simql \Delta'$.
Then 
$$
\depth_x\o_X(-D)\geq \min\{3,\codim_Xx\}.
$$
\item  Let $Z\subset X$ be any closed,  reduced  subscheme that is a union of
lc centers of $(X,\Delta)$. Then
$$
\depth_x\o_X(-Z)\geq \min\{3,1+\codim_Zx\}.
$$
\end{enumerate}
\end{thm}

 In contrast with 
 (\ref{d.num.eq.to.Delta.CM.thm}), my proof does not work in the
complex analytic case; see (\ref{anal.ref.rem}).

\begin{say}[Applications and examples]\label{d.num.eq.to.Delta.CM.appl}{\ }

(\ref{d.num.eq.to.Delta.CM.appl}.1)
In (\ref{d.num.eq.to.Delta.CM.thm})  
we can take $\Delta'=0$. Then  $D$ can be any $\q$-Cartier divisor,
reproving  \cite[5.25]{km-book}. The  $D=\Delta'\leq \Delta$ case
recovers 
\cite[4.13]{fujinobook}.

(\ref{d.num.eq.to.Delta.CM.appl}.2)
The $D=0$ case of (\ref{main.slc.thm}.1)
is a theorem of  \cite{ale-lim, fujinobook} which says that
if  $(X,\Delta)$ is lc  and $x\in X$  is not a log canonical center
then $\depth_x\o_X\geq \min\{3,\codim_Xx\}$.
This can fail if $x$ is a log canonical center, for instance
when $x\in X$ is a cone over an Abelian variety of dimension $\geq 2$.

(\ref{d.num.eq.to.Delta.CM.appl}.3) If  $(X,\Delta)$ is slc 
then $K_X+\Delta$ is $\q$-Cartier, hence
$-K_X\simql \Delta$. 
Then $\o_X\bigl(-(-K_X)\bigr)\cong \omega_X$.
Thus if  $x\in X$  is not a log canonical center
then   $\depth_x\omega_X\geq \min\{3,\codim_Xx\}$.
(Note that while $\o_X$ is CM iff $\omega_X$ is CM, it can happen
that $\o_X$ is $S_3$ but $\omega_X$ is not. Thus 
 (\ref{d.num.eq.to.Delta.CM.appl}.3) does not seem to be a formal
consequence of (\ref{d.num.eq.to.Delta.CM.appl}.2).)

Let $f:(X,\Delta)\to C$ be a semi log canonical morphism
to a smooth curve $C$ (cf.\ \cite[7.1]{km-book}) and
 $X_c$ the fiber over a closed point.
None of the lc centers are contained in $X_c$,
thus if $x\in X_c$ has codimension $\geq 2$
then $\depth_x\omega_{X/C}\geq 3$. Therefore, 
the restriction of $\omega_{X/C}$
to $X_c$ is $S_2$, hence it is isomorphic to $\omega_{X_c}$.
More generally,  $\omega_{X/C}$ commutes with arbitrary base change.
(When the general fiber is klt, this follows from 
\cite{elkik};
for projective morphisms a proof is given in \cite{k-k-db},
but the general case has not been known earlier.
As far as I know, the complex analytic case is still unproved.) 

(\ref{d.num.eq.to.Delta.CM.appl}.4)  Assume that $(X,\Delta)$ is slc.
For any $n\geq 1$, write
$$
-nK_X-\rdown{n\Delta}\simq -n\bigl(K_X+\Delta)+
\bigl(n\Delta-\rdown{n\Delta}\bigr)\simql
n\Delta-\rdown{n\Delta}.
$$
Assume now that
$\Delta=\sum \bigl(1-\frac1{m_i}\bigr)D_i$ with $m_i\in \n\cup\{\infty\}$.
Then $n\Delta-\rdown{n\Delta}=\sum_i \frac{c_i}{m_i}D_i$
for some $c_i\in\n$ 
where $0\leq c_i<m_i$ for every $i$. Thus $c_i\leq m_i-1$ for every $i$,
that is, 
  $n\Delta-\rdown{n\Delta}\leq \Delta$. Thus, if
 $x\in X$  is not a log canonical center
then 
$$
\depth_x\o_X\bigl(nK_X+\rdown{n\Delta}\bigr)\geq \min\{3,\codim_Xx\}.
$$
In particular, if $f:(X,\Delta)\to C$ is a proper slc morphism
to a smooth curve and $X_c$ is any fiber then
the restriction of $\o_X\bigl(nK_X+\rdown{n\Delta}\bigr)$
to any fiber is $S_2$ and hence the natural map
$$
\o_X\bigl(nK_X+\rdown{n\Delta}\bigr)|_{X_c}\to
\o_{X_c}\bigl(nK_{X_c}+\rdown{n\Delta_c}\bigr)
\qtq{is an isomorphism.}
$$
This implies that  the
Hilbert function of the fibers 
$$
\chi\bigl(X_c,\o_{X_c}\bigl(nK_{X_c}+\rdown{n\Delta_c}\bigr)\bigr)
$$
is deformation invariant.
(Note that, because of the rounding down, the Hilbert function
is not a polynomial, rather a polynomial whose coefficients are
periodic functions of $n$. The period divides the index of
$(X,\Delta)$, that is, the
smallest $n_0\in \n$ such that $n_0\Delta$ is a $\z$-divisor and
$n_0\bigl(K_X+\Delta\bigr)$ is Cartier.)

(\ref{d.num.eq.to.Delta.CM.appl}.5)
It is also  worthwhile to note that while the assumptions of
Theorem \ref{d.num.eq.to.Delta.CM.thm} depend only on the
$\q$-linear equivalence class of $D$, 
being CM is not preserved by
$\q$-linear equivalence in  general. For instance, let $X$ be a cone
over an Abelian variety $A$ of dimension $\geq 2$.
Let $D_A$ be a $\z$-divisor on $A$ such that $mD_A\sim 0$
for some $m>1$ but $D_A\not\sim 0$. Let $D_X$ be the cone over $D_A$.
Then $D_X\simql 0$, $\o_X(D_X)$ is CM but $\o_X$ is not CM.

These assertions  follow from the next easy characterization of CM 
divisorial sheaves on cones:

(\ref{d.num.eq.to.Delta.CM.appl}.6)
 {\it Claim.} Let $Y\subset \p^n$ be projectively normal, 
$H$ the hyperplane class on $Y$  and $D$ a 
Cartier divisor on $Y$. Let $X \subset \a^{n+1}$ be the cone over $Y$
with vertex $v$ 
and $D_X$ the cone over $D$. Then
$$
H^i_v\bigl(X, \o_X(D_X)\bigr)=
\sum_{m\in\z} H^{i-1}\bigl(Y, \o_Y(D+mH)\bigr)
\qtq{for $i\geq 2.$}
$$
In particular,  $\o_X(D_X)$ is CM iff
$$
H^i\bigl(Y, \o_Y(D+mH)\bigr)=0\quad \forall\ m\in \z, \ \forall\  0<i<\dim Y.
\qed
$$

(\ref{d.num.eq.to.Delta.CM.appl}.7)
 Consider the quadric cone $X:=(x_1x_2=x_3x_4)\subset \a^4$
with vertex $v=(0,0,0,0)$.
It is the cone over the quadric surface $Q\cong\p^1\times \p^1 \subset \p^3$.
It contains two families of planes with typical members
$A:=(x_1=x_3=0)$ and $B:=(x_1=x_4=0)$.
By (\ref{d.num.eq.to.Delta.CM.appl}.6)
$$
\begin{array}{l}
H^2_v\bigl(X, \o_X(aA+bB)\bigr)=
\sum_{m\in\z} H^1\bigl(Q, \o_Q(a+m,b+m)\bigr)\\[1ex]
\qquad\qquad=
\sum_{0\leq m\leq |b-a|-2} 
H^0\bigl(\p^1, \o_{\p^1}(m)\bigr)\otimes
H^0\bigl(\p^1, \o_{\p^1}(|b-a|-2-m)\bigr).
\end{array}
$$
Thus we see that $\o_X(aA+bB)$ is CM only if $|b-a|<2$.

(\ref{d.num.eq.to.Delta.CM.appl}.8)
As another application of (\ref{d.num.eq.to.Delta.CM.thm}),
 assume that  $\left(X, \sum a_iD_i\right)$ is
dlt and $1-\frac1{n}\leq a_i\leq 1$ for every $i$ for some $n\in \n$. 
Then, for every $m$, 
$$
\begin{array}{rll}
m\bigl(K_X+\tsum_i D_i\bigr)&=&  \tsum_i m(1-a_i)D_i+ 
m\bigl(K_X+\tsum_i a_iD_i\bigr)\\
&\simql &
  \tsum_i m(1-a_i)D_i.
\end{array}
$$
If $1\leq m\leq n-1$ then 
 $0\leq m(1-a_i)\leq a_i$, thus 
by (\ref{d.num.eq.to.Delta.CM.thm}) and by Serre duality
we conclude that
$$
\omega_X^{[-m]}\left(-m\tsum D_i\right)
\qtq{and}
\omega_X^{[m+1]}\left(m\tsum D_i\right)
\qtq{are CM for $1\leq m\leq n-1$.}
\eqno{(\ref{d.num.eq.to.Delta.CM.appl}.8.1)}
$$
If, in addition,   $K_X$ is $\q$-Cartier, then
$m\tsum_i D_i\simql \tsum_im(1-a_i)D_i$, hence 
$$
\o_X\left(-m\tsum D_i\right)\qtq{is CM for $1\leq m\leq n-1$.}
\eqno{(\ref{d.num.eq.to.Delta.CM.appl}.8.2)}
$$
Results like these are quite fragile.
As an example, let $X\subset \a^4$ be the quadric cone
with the 2 families of planes $|A|$ and $|B|$. 
Then 
$$
\bigl(X, A_1+\tfrac12(B_1+B_2)\bigr)\qtq{and}
\bigl(X, \tfrac9{10}(A_1+A_2)+\tfrac6{10}(B_1+B_2+B_3)\bigr)
$$
are both dlt, giving that
$$
\o_X\left(-A_1-B_1-B_2\right)\qtq{and}
\o_X\left(-A_1-A_2-B_1-B_2-B_3\right)\qtq{are CM.}
$$
Note, however, that the sheaves
$$
\o_X\left(-B_1-B_2\right), \o_X\left(-B_1-B_2-B_3\right)\qtq{and}
\o_X\left(-A_1-B_1-B_2-B_3\right)
$$
are not CM by (\ref{d.num.eq.to.Delta.CM.appl}.7).

(\ref{d.num.eq.to.Delta.CM.appl}.9)
The following example shows that
(\ref{d.num.eq.to.Delta.CM.appl}.8.2) fails in general if
$K_X$ is not $\q$-Cartier. 

Let $Q\subset\a^4$ be the affine quadric $(xy=zt)$.
Let $B_1=(x=z=0)$ and $B_2=(y=t=0)$ be 2 planes in  the same family
of planes on $Y$.
For some $c_i$  (to be specified later),
consider the divisor $c_1B_1+c_2B_2$.
(Note that  $K_Y+c_1B_1+c_2B_2$ is $\q$-Cartier only if
$c_1+c_2=0$, so one of the $c_i$  will have to be
negative, but let us not worry about it for now.)

Consider the
group action $\tau: (x,y,z,t)\mapsto (\epsilon x,y,\epsilon z,t)$
where $\epsilon $ is a primitive $n$th root of unity.
This generates an action of $\mu_n$; let
$\pi:Y\to X_{n}=Y/\mu_n$ be the corresponding quotient.
Both of the $B_i$ are $\tau$-invariant; set $B'_i=\pi(B_i)$
(with reduced structure).

Note that the fixed point set of $\tau$ is $B_1$ and 
$\pi$ ramifies along $B_1$ with ramification index $n$.
Thus $K_Y=\pi^*\bigl(K_{X_{n}}+(1-\tfrac1{n})B'_1\bigr)$ hence
$$
K_Y+c_1B_1+c_2B_2=\pi^*\Bigl(K_{X_{n}}+
\bigl(1-\tfrac1{n}+\tfrac{c_1}{n}\bigr)B'_1+ c_2B'_2\Bigr).
$$
Now we see that even if $c_1<0$, the
coefficient of $B'_1$ could be positive.
In particular one computes that
$$
K_Y-\bigl(1-\tfrac2{n+1}\bigr)B_1+\bigl(1-\tfrac2{n+1}\bigr)B_2=
\pi^*\Bigl(K_{X_{n}}+
\bigl(1-\tfrac2{n+1}\bigr)\bigl(B'_1+B'_2\bigr)\Bigr).
$$
Thus $\Bigl({X_{n}},
\bigl(1-\tfrac2{n+1}\bigr)\bigl(B'_1+B'_2\bigr)\Bigr)$
is klt but 
 $B'_1+B'_2$ consists of 2  normal surfaces
intersecting at a point, hence it is not $S_2$.
Therefore, $\o_{X_n}\bigl(-B'_1-B'_2\bigr)$ has only depth 2 at the origin.

(\ref{d.num.eq.to.Delta.CM.appl}.10)
With $D$ as in  (\ref{d.num.eq.to.Delta.CM.thm}), $\o_X(-D)$ is CM. 
 What about $\o_X(D)$?

On a proper variety, sheaves of the form
$\o_X(D)$ are quite different from ideal sheaves, but being CM
is a local condition. If $X$ is affine, then
there is always a reduced divisor $D'$ such that
 $\o_X(D)\cong \o_X(-D')$. Despite this, (\ref{d.num.eq.to.Delta.CM.thm})
does not hold for $\o_X(D)$.

To see such an example, let $S:=\p^1\times \p^1$ 
and $C:=\p^1\times \{(0{:}0)\}$ a line on $S$. Embed $S$ into $\p^5$ by
$|-(K_S+C)|$  and let $X\subset \a^6$ be the affine cone over $S$
and $D\subset X$ the cone over $C$. Then $(X,D)$ is  canonical
and so $\o_X(-D)$ is CM.

Let  $F_i\subset X$ be cones over lines 
of the form $ \{p_i\}\times\p^1$. Then $D+F_1+F_2\sim 0$, hence
 $\o_X(D)$ is isomorphic to
$\o_X(-F_1-F_2)$.
Since $F_1+F_2$ is not $S_2$,   $\o_X(D)$ is not CM.
\end{say}

The proof of Theorem \ref{d.num.eq.to.Delta.CM.thm}
uses the method of  two spectral sequences
introduced in \cite[5.22]{km-book} in the global case
and in \cite{fujinobook} in the local case.

\begin{say}[The method of  two spectral sequences]\label{d.two.ss.say}

Let $f:Y\to X$ be a proper morphism, 
$V\subset X$ a closed subscheme and
$W:=f^{-1}V\subset Y$. For any coherent sheaf $F$  on $Y$
there is a  Leray spectral sequence
$$
\begin{array}{lcl}
H^i_V\bigl(X, R^jf_*F\bigr)& \Rightarrow & H^{i+j}_W(Y,F),
\end{array}
\eqno{(\ref{d.two.ss.say}.1)}
$$
where $H^i_V$ denotes cohomology with supports in $V$;
see \cite[Chap.1]{sga2}. 
In particular, if $R^jf_*F=0$ for every $j>0$ then the
spectral sequence degenerates and we get isomorphisms
$H^i_V\bigl(X, f_*F\bigr)\cong H^{i}_W(Y,F)$.

Given a map of sheaves $F\to F'$ we get, for each $i$,  a commutative diagram
$$
\begin{array}{ccc}
 H^i_V\bigl(X, f_*F'\bigr)  & \stackrel{\alpha'_i}{\to} & 
H^i_W\bigl(Y, F'\bigr)\\
\uparrow & & \uparrow \\
 H^i_V\bigl(X, f_*F\bigr)  & \stackrel{\alpha_i}{\to} & H^i_W\bigl(Y, F\bigr).
\end{array}
\eqno{(\ref{d.two.ss.say}.2)}
$$
The following simple observation will be a key ingredient in the proof of
(\ref{d.num.eq.to.Delta.CM.thm}).
\end{say}

\begin{claim} \label{d.two.ss.say.3} With the above notation, assume that
for some $i$, 
\begin{enumerate}
\item $f_*F=f_*F'$, 
\item $H^i_W\bigl(Y, F\bigr)=0$ and
\item $\alpha'_i$ is an isomorphism.
\end{enumerate}
Then $H^i_V\bigl(X, f_*F\bigr)=0$. \qed
\end{claim}

For the cases $i=1,2$, somewhat weaker hypotheses suffice. 
First note that
$H^1_V\bigl(X, f_*F\bigr)\into H^1_W\bigl(Y, F\bigr)$ is injective,
thus, if $i=1$, then (\ref{d.two.ss.say.3}.2) alone
yields $H^1_V\bigl(X, f_*F\bigr)=0$.

The $i=2$ case is more interesting.
If $H^2_W\bigl(Y, F\bigr)=0$ then $\alpha_2=0$. 
Together with  (\ref{d.two.ss.say.3}.1) this
implies that $\alpha'_2=0$. On the other hand, $\alpha'_2$
sits in the exact sequence
$$
  H^0_V\bigl(X, R^1f_*F'\bigr)\to H^2_V\bigl(X, f_*F'\bigr)
\stackrel{\alpha'_2}{\longrightarrow} H^2_W\bigl(Y, F'\bigr),
$$
hence we get the following, first used in \cite{ale-lim}.

\begin{claim}\label{d.two.ss.say.4} With the above notation, assume that 
\begin{enumerate}
\item $f_*F=f_*F'$, 
\item $H^2_W\bigl(Y, F\bigr)=0$ and
\item $H^0_V\bigl(X, R^1f_*F'\bigr)=0$.
\end{enumerate}
Then $H^2_V\bigl(X, f_*F\bigr)=0$. \qed
\end{claim} 

\begin{say}[A special case]\label{d.two.ss.II.say}
As a warm up we prove,
following the methods of \cite[5.22]{km-book} and
\cite[4.2.1.App]{fujinobook},
 that if   $(X,\Delta)$ is klt then
$X$ is CM and has rational singularities.

We need to prove that
$H^i_x(X,F)=0$ for $i<\codim_Xx$ for every point $x\in X$.
We can localize at $x$; thus from now on assume that $x$ is a closed point.
Let $f:Y\to X$ be a resolution such that
$W:=f^{-1}(x)$ is a divisor.
Choose $F:=\o_Y$ and $F':=\o_Y(B)$ where $B$ is an effective, $f$-exceptional
divisor to be specified later.
Then $f_*\o_Y(B)=f_*\o_Y=\o_X$ 
hence (\ref{d.two.ss.say.3}.1) holds
and  we have
a commutative diagram
$$
\begin{array}{ccc}
 H^i_x\bigl(X, \o_X\bigr)  & \stackrel{\alpha'_i}{\to} & 
H^i_W\bigl(Y, \o_Y(B)\bigr)\\
|| & & \uparrow \\
 H^i_x\bigl(X, \o_X\bigr)  & \to & 
H^i_W\bigl(Y,  \o_Y\bigr).
\end{array}
\eqno{(\ref{d.two.ss.II.say}.1)}
$$
By (\ref{duality.cor}), $H^i_W\bigl(Y,  \o_Y\bigr)=0$ for $i<\dim X$,
hence (\ref{d.two.ss.say.3}.2) also holds.

In order to prove  (\ref{d.two.ss.say.3}.3), we finally use
that  $(X,\Delta)$ is klt.
By definition, this means than we can choose $f:Y\to X$
such that
$$
 K_{Y} \simq f^*(K_X+\Delta) +B-A
\eqno{(\ref{d.two.ss.II.say}.2)}
$$
 where
$B$ is  an effective, $f$-exceptional, $\z$-divisor,
$A$ is a simple normal crossing divisor
   and $\rdown{A}=0$.
Then  
$B\simq K_Y+A-f^*(K_X+\Delta)$, hence we conclude from
(\ref{relkod.klt.thm}) that 
$R^if_*\o_Y(B)=0$  for $i>0$.
Therefore the  spectral sequence for $\o_Y(B)$ degenerates and
$H^i_x\bigl(X, \o_X\bigr)   \cong  
H^i_W\bigl(Y, \o_Y(B)\bigr)$ for every $i$.

Thus, for $i<\dim X$,  the 
commutative diagram (\ref{d.two.ss.II.say}.1) becomes
$$
\begin{array}{ccc}
 H^i_x\bigl(X, \o_X\bigr)  & \cong & 
H^i_W\bigl(Y, \o_Y(B)\bigr)\\
|| & & \uparrow \\
 H^i_x\bigl(X, \o_X\bigr)  & \to & \hphantom{.}0.
\end{array}
\eqno{(\ref{d.two.ss.II.say}.3)}
$$
This implies that $ H^i_x\bigl(X, \o_X\bigr)=0$
for $i<\dim X$, hence $X$ is CM. 

Next we prove that $R^jf_*\o_Y=0$ for $j>0$.
By induction on the dimension and localization, we may
assume that $ \supp R^jf_*\o_Y\subset \{x\}$  for $j>0$.
Then 
$$
H^i_x\bigl(X, R^jf_*\o_Y\bigr)=0
\qtq{unless $i=0$ or $(i,j)=(n,0)$.}
$$
Since $H^i_W\bigl(Y,  \o_Y\bigr)=0$ for $i<\dim X$,
we conclude that $R^jf_*\o_Y=0$ for $0<j<n-1$ and 
we have an exact sequence
$$
0\to R^{n-1}f_*\o_Y\to H^n_x\bigl(X, \o_X\bigr)
\stackrel{\alpha_n}{\to} H^n_W\bigl(Y,  \o_Y\bigr).
\eqno{(\ref{d.two.ss.II.say}.4)}
$$
Note that $\alpha_n$ also sits in the diagram
$$
\begin{array}{ccc}
 H^n_x\bigl(X, \o_X\bigr)  & \cong & 
H^n_W\bigl(Y, \o_Y(B)\bigr)\\
|| & & \uparrow \\
 H^n_x\bigl(X, \o_X\bigr)  & \stackrel{\alpha_n}{\to}
 & \hphantom{.}H^n_W\bigl(Y,  \o_Y\bigr)
\end{array}
\eqno{(\ref{d.two.ss.II.say}.5)}
$$
which implies that $\alpha_n$ is injective.
Thus $R^{n-1}f_*\o_Y=0$ as required. \qed
\end{say}

The proof of the general case is quite similar.

\begin{say}[Proof of  Theorem \ref{d.num.eq.to.Delta.CM.thm}]
\label{d.num.eq.to.Delta.CM.pf}{\ }

We may assume that $X$ is affine
and $x\in X$ is a closed point. Write $\Delta=\Delta'+\Delta''$.
As in \cite[2.43]{km-book}, there are effective $\q$-divisors
$\Delta'_1\simql \Delta'$ and $\Delta''_1\simql \Delta''$ such that 
$$
\bigl(X, (1-\epsilon)\Delta+\epsilon(\Delta'_1+\Delta''_1)\bigr)
\qtq{is klt.}
$$
Furthermore, $D\simql \Delta'\simql (1-\epsilon)\Delta'+\epsilon\Delta'_1$.
Thus we may assume that $(X,\Delta)$ is in fact klt.

We need to prove that for every $x\in X$,
$$
H^i_x\bigl(X,\o_X(-D)\bigr)=0\qtq{for $i<\codim_Xx$.}
$$
We can localize at $x$, hence we may assume that
$x$ is a closed point.

Choose  a log resolution  $f:Y\to X$ of  $(X,\Delta)$
such that $W:=f^{-1}(x)$ is a divisor  and write
$$
f^*\bigl(D-\Delta'\bigr)=f^{-1}_*D-f^{-1}_*\Delta'-F,
\eqno{(\ref{d.num.eq.to.Delta.CM.pf}.1)}
$$
where $F$ is $f$-exceptional. Set
$D_Y:=f^{-1}_*D-\rdown{F}$ and note that
$$
D_Y=f^{-1}_*\Delta'+\{F\}+f^*\bigl(D-\Delta'\bigr).
$$
Thus, from (\ref{d.S2.pushfwd.lem}), we obtain that 
 if $B_Y$ is effective and $f$-exceptional, then
$f_*\o_Y\bigl(B_Y-D_Y\bigr)=\o_X(-D)$.
Thus (\ref{d.two.ss.say.3}.1) holds and we have a commutative diagram
$$
\begin{array}{ccc}
 H^i_x\bigl(X, \o_X(-D)\bigr)  & \stackrel{\alpha'_i}{\to} & 
H^i_W\bigl(Y, \o_Y(B_Y-D_Y)\bigr)\\
|| & & \uparrow \\
 H^i_x\bigl(X, \o_X(-D)\bigr)  & \to & 
H^i_W\bigl(Y,  \o_Y(-D_Y)\bigr)
\end{array}
\eqno{(\ref{d.num.eq.to.Delta.CM.pf}.2)}
$$

By (\ref{d.num.eq.to.Delta.CM.pf}.1), 
$D_Y\simq  \bigl(\mbox{$f$-nef})+\bigl(f^{-1}_*\Delta'+\{F\}\bigr)$,
thus, by (\ref{duality.cor}), 
$H^i_W\bigl(Y, \o_Y(-D_Y)\bigr)=0$ for $i< \dim X$
and so (\ref{d.two.ss.say.3}.2)  is satisfied.

Finally we need $\alpha'_i$ to be an isomorphism.
By the klt assumption, 
$$
 K_{Y} +f^{-1}_*\Delta\simq f^*(K_X+\Delta) +B-A
\eqno{(\ref{d.num.eq.to.Delta.CM.pf}.3)}
$$
 where
$A, B$ are effective, $f$-exceptional, $B$ is a $\z$-divisor,
  and $\rdown{A}=0$.
 Write 
$B-A+\{F\}=:B_Y-A_Y $ where 
$A_Y, B_Y$ are effective, $f$-exceptional, $B_Y$ is a $\z$-divisor,
  and $\rdown{A_Y}=0$. Note that
$$
\begin{array}{rcl}
B_Y-D_Y&\simq &B-A+\{F\}+A_Y-f^{-1}_*\Delta'-\{F\}-f^*\bigl(D-\Delta'\bigr)\\
&\simq& K_{Y} +f^{-1}_*\Delta +A_Y-f^{-1}_*\Delta'- 
f^*\bigl(K_X+\Delta+D-\Delta'\bigr)\\
&\simq &K_{Y} +f^{-1}_*\Delta'' +A_Y- 
f^*\bigl(K_X+\Delta+D-\Delta'\bigr).
\end{array}
\eqno{(\ref{d.num.eq.to.Delta.CM.pf}.4)}
$$
Thus $R^if_*\o_Y\bigl(B_Y-D_Y\bigr)=0$ for $i>0$ 
by (\ref{relkod.klt.thm}) and so
$\alpha'_i$ is  an isomorphism.

These imply that $ H^i_x\bigl(X, \o_X(-D)\bigr)=0$ for $i<\dim X$
hence $\o_X(-D)$ is CM. \qed
\end{say}

\begin{say}[Proof of  Theorem \ref{main.slc.thm}]
\label{main.slc.thm.pf}{\ }

We start with (\ref{main.slc.thm}.1) and
first consider the case when $X$ is normal, that is, when $(X,\Delta)$ is lc.
There are a few places where we have to modify the previous proof
(\ref{d.num.eq.to.Delta.CM.pf}).

We may assume that $\rdown{\Delta}=0$.
 (If $\Delta=\sum d_iD_i$  then we  can replace
$\Delta$ by $\sum (d_i/2)(D_{1i}+D_{2i})$ where
$D_{1i}, D_{2i}$ are general members of the linear system
$|D_i|$, cf.              \cite[2.33]{km-book}.)
Let $f:Y\to X$ be a log resolution of $(X,\Delta)$ 
and write
$$
 K_{Y} +f^{-1}_*\Delta\simq f^*(K_X+\Delta) +B-A-E,
\eqno{(\ref{d.num.eq.to.Delta.CM.pf}.5)}
$$
 where
$A, B,E$ are effective, $f$-exceptional, $B,E$ are $\z$-divisors,
$E$ is reduced
  and $\rdown{A}=0$.

Pick $D_Y$ as in (\ref{d.num.eq.to.Delta.CM.pf}.1). By (\ref{d.S2.pushfwd.lem}),
 if $B_Y$ is effective and $f$-exceptional, then
$f_*\o_Y\bigl(B_Y-D_Y\bigr)=\o_X(-D)$
and $H^i_W\bigl(Y, \o_Y(-D_Y)\bigr)=0$ for $i<\dim X$ by
(\ref{relkod.klt.thm}).
It remains to check (\ref{d.two.ss.say.4}.3), that is,
the vanishing of
$H^0_x\bigl(X,  R^1f_*\o_Y\bigl(B_Y-D_Y\bigr)\bigr)$.

To this end choose  $B_Y,A_Y,E_Y$ such that 
$$
B_Y-A_Y -E_Y=B-A-E+\{F\}.
$$
  where
$A_Y, B_Y, E_Y$ are effective, $f$-exceptional, $B_Y, E_Y$ are $\z$-divisors,
$E_Y$ is reduced, 
  $\rdown{A_Y}=0$ and $F$ as in (\ref{d.num.eq.to.Delta.CM.pf}.1).
 With this choice, as in
(\ref{d.num.eq.to.Delta.CM.pf}.4),
$$
B_Y-D_Y\simq K_{Y} +f^{-1}_*\Delta'' +A_Y+E_Y- 
f^*\bigl(K_X+\Delta+D-\Delta'\bigr).
$$
By  assumption $x$ is not a log canonical center, hence by
(\ref{amb-fuj.thm}) $x$ is not an  associated prime
of $R^if_*\o_Y(B_Y-D_Y)$. Thus 
 $H^0_x\bigl(X,  R^1f_*\o_Y(B_Y-D_Y)\bigr)=0$
 and hence
$H^2_x\big(X, \o_X(-D)\bigr) =0$.

The above proof should work without changes if 
$X$ is not normal, that is, when $(X,\Delta)$ is slc, but 
(\ref{amb-fuj.thm}) is not stated for semi-resolutions in the references.
We go around this problem as follows. 

By (\ref{double.cover.trick}) there is a double cover
$\pi: (\tilde X,\tilde \Delta)\to (X,\Delta)$, \'etale in codimension 1
such that every irreducible component of $(\tilde X,\tilde \Delta)$
is smooth in codimension 1. Set $\tilde D:=\pi^{-1}(D)$.
Then $\o_X(-D)$ is a direct summand of 
$\pi_*\o_{\tilde X}(-\tilde D)$, hence it is enough to prove
the depth bounds for $\o_{\tilde X}(-\tilde D)$.

As in \cite[20]{k-semires} we can construct a semi-resolution
 $\tilde f:(\tilde Y, \tilde \Delta_Y)\to (\tilde X,\tilde \Delta)$  
such that $(\tilde Y, \tilde \Delta_Y)$ is an embedded simple normal 
crossing pair, as required in (\ref{amb-fuj.thm}). 
The rest of the proof works as before.

Next consider (\ref{main.slc.thm}.2).
The only interesting case is when $\codim_Zx\geq 2$.
Since $x$ is not a lc center, this implies that $\codim_{Z_i}x\geq 2$
for every irreducible component $Z_i\subset Z$. 
(It is, however, possible that $x$ has codimension 1 in some
lc center that is contained in $Z$.)
By localizing, we may assume that $x$ is a closed point.

If $X$ is normal, let $f:Y\to X$ be a  log resolution.
If $X$ is not normal, as before, by first passing to a double cover
(\ref{double.cover.trick}) we may assume that there is a
 semi log resolution  $f_1:Y_1\to X$ to which 
 (\ref{amb-fuj.thm})  applies.
We may also assume that 
every irreducible component of $Z$ is the image of
some divisor $E_j\subset Y$ with discrepancy $-1$.
We have to be more careful if $Z$ contains one of the
codimension 1  components of $\sing X$. In this case,
we have a divisor 
$E_j\subset \sing Y_1$ mapping to $Z$. We blow up $E_j$  to get
$f:Y\to X$ and replace $E_j$ with both of the irreducible components of
its preimage on $Y$.   
Set $E=\sum_j E_j$. By  (\ref{pushfwd.easy.cond}),
$$
\o_X(-Z)=f_*\bigl(\o_Y(-E)\bigr)=f_*\bigl(\o_Y(B-E)\bigr)
$$
for any effective $f$-exceptional divisor $B$ whose support
 does not contain any of the $E_j$.
As usual, write
$$
 K_{Y} \simq f^*(K_X+\Delta) +B-E-E'-A.
$$
By (\ref{duality.prop}), $H^i_W\bigl(Y,\o_Y(-E)\bigr)=0$
is dual to $\bigl(R^{n-i}f_*\omega_Y(E)\bigr)_x$.
We have assumed that $\dim f(E_j)\geq 2$ for every $j$ and
 $\dim f(E_j\cap E_k)\geq 1$ for every $j\neq k$ since
$x$ is not a lc center. Thus, by (\ref{relkod.lc.cor}) and 
(\ref{duality.prop}),
 $H^i_W\bigl(Y,\o_Y(-E)\bigr)=0$ for $i\leq 2$
and so   (\ref{d.two.ss.say.4}.2) holds.

From 
$B-E\simq K_Y+E'+A-f^*(K_X+\Delta)$
and (\ref{amb-fuj.thm}) we conclude that every associated prime
of $R^if_*\o_Y(B-E)$ is  an lc center of $(X,\Delta)$.
Since $x$ is not an  lc center, this implies that
 $H^0_x\bigl(X, R^1f_*\o_Y(B-E)\bigr)=0$,
giving (\ref{d.two.ss.say.4}.3).
Therefore
$H^2_x\big(X, \o_X(-Z)\bigr) =0$ and
so  $\o_X(-Z)$ has depth $\geq 3$ at $x$.
\qed
\end{say}

\subsection*{Conditions for (\ref{d.two.ss.say.3}.1)}{\ }

\begin{say}\label{pushfwd.easy.cond}
 The assumption  (\ref{d.two.ss.say.3}.1)
is easy to satisfy in many cases. 
Let $f:Y\to X$ be a proper, birational morphism to a normal scheme $X$.
For any closed subscheme
$Z_Y\subset Y$, $f_*\o_Y(-Z_Y)=\o_X\bigl(-f(Z_Y)\bigr)$
where $f(Z_Y)$ is the scheme theoretic image.
If $Z_Y$ is reduced, then  $Z_X:=f(Z_Y)\subset X$
is also reduced. Thus if 
  $B$ is $f$-exceptional and  $f(B)$ does not contain any of the
irreducible components of $Z_X$, then 
$\o_X(-Z_X)=f_*\o_Y(-Z_Y)=f_*\o_Y(B-Z_Y)$.

The last equality holds even if $X$ is not normal.
\end{say}

Another easy  case is the following (cf.\ \cite{fujita85}).

\begin{lem}\label{d.S2.pushfwd.lem} 
Let $f:Y\to X$ be a proper, birational morphism.
Let $D$ be a $\z$-divisor on $X$ and assume that
$D\sim_{\q,f}D_h+D_v$ where $D_v$ is effective, $f$-exceptional
and $D_h$ is effective without exceptional components.
Let $B$ be an  effective, $f$-exceptional divisor
whose support does not contain any of the irreducible
components of  $\rdown{D_v}$.
Assume that
\begin{enumerate}
\item either  $X$ and $Y$ are normal,
\item or $X$ and $Y$ are $S_2$, $f$ is an isomorphism outside a
codimension 2 subscheme of $X$ and $Y$ is normal at the generic
point of every exceptional divisor.
\end{enumerate}
Then 
$$
\o_X\bigl(-f_*D\bigr)=f_*\o_Y(-D)=f_*\o_Y(B-D).
$$
\end{lem}

Proof. Fix a  section  $s\in H^0\bigl(X, f_*\o_Y(B-D)\bigr)$ and
choose $B$ as small as possible.

By assumption, there is a $\q$-Cartier $\q$-divisor $M$ on $X$ such that
$D+f^*M=D_h+D_v$. Choose $n\in\n$ such that $nD_h, nD_v, f^*(nM)$ are
 all $\z$-divisors.
Then 
$$
\begin{array}{rcl}
\bigl(f^*s\bigr)^n&\in&  H^0\bigl(Y, \o_Y(nB-nD)\bigr)\\
&=& H^0\bigl(Y, \o_Y(nB-nD_h-nD_v+ f^*(nM)\bigr)\\
&\subset &
H^0\bigl(Y, \o_Y(nB-nD_v+ f^*(nM)\bigr).
\end{array}
$$
As noted in (\ref{pushfwd.easy.cond}), 
adding effective exceptional divisors to a pull-back never creates
new sections. By assumption,  every irreducible component of $B$
appears in $nB-nD_v$ with positive coefficient. 
 Thus $\bigl(f^*s\bigr)^n$ vanishes along every  irreducible component of $B$,
contradicting the minimality of $B$. Thus $B=0$ and so
$s$ is a   section of $f_*\o_Y(-D)$.\qed

\subsection*{Conditions for (\ref{d.two.ss.say.3}.2)}{\ }

The following is the relative version
of the  Kawamata--Viehweg vanishing theorem
\cite{kaw, vie}.

\begin{thm}\label{relkod.klt.thm}
 Let $f:Y\to X$ be a projective, birational  morphism, $Y$ smooth
and $\Delta$ an effective simple normal crossing divisor on $Y$
such that $\rdown{\Delta}=0$.
Let $M$ be a line bundle on $Y$ and assume that
$M\sim_{\q,f}K_Y+(\mbox{$f$-nef})+\Delta$. Then
$R^if_*M=0$ for $i>0$. \qed
\end{thm}

If $f$ is not birational, $Y$ is not smooth or if $\rdown{\Delta}\neq 0$, then
vanishing fails in general. 
There are, however, some easy consequences
that can be read off from induction on the number of
irreducible components and by writing down the obvious exact sequences
$$
0\to  M\to  M(H) \to  M(H)|_H\to 0
\qtq{and}
0\to  M(-D)\to  M \to  M|_D\to 0
$$
where $H$ is a smooth  sufficiently ample divisor on $Y$ and 
$D\subset \rdown{\Delta}$ is a smooth divisor.

Let $W$ be a smooth variety and $\sum_{i\in I} E_i$ a snc divisor on $W$.
Write $I=I_V\cup I_D$ as a disjoint union. Set $Y:=\sum_{i\in I_V} E_i$
as a subscheme and $D_Y:=\sum_{i\in I_D} E_i|_Y$ as a divisor on $Y$.
We call $(Y, D_Y)$ an {\it embedded snc pair.}
Anything isomorphic to such a pair is called an
{\it embeddable snc pair.} A pair is called an {\it  snc pair}
if it is locally an embeddable snc pair.

\begin{cor}\label{relkod.klt.nonbir.thm}
 Let $(Y, \sum_iD_i)$  be an snc pair
and  $f:Y\to X$ be a projective morphism.
Let $M$ be a line bundle on $Y$ and assume that
$$
M\sim_{\q,f}K_Y+(\mbox{$f$-nef})+\tsum_i a_i D_i\qtq{where $0\leq a_i\leq 1$.}
$$
Then $R^if_*M=0$ for $i>\dim Y-\dim X$. \qed
\end{cor}

\begin{cor}\label{relkod.lc.cor}
Notation and assumptions as in (\ref{relkod.klt.nonbir.thm}).
Set $\Delta=\sum_i a_i D_i$ and $n=\dim Y$. Then
\begin{enumerate}
\item $R^nf_*M=0$,
\item $R^{n-1}f_*M=0$ unless there is a
divisor $B\subset \sing Y\cup\rdown{\Delta}$ such that $\dim f(B)=0$ and
\item $R^{n-2}f_*M=0$ unless 
there are
divisors  $B_1, B_2\subset \sing Y\cup\rdown{\Delta}$ such that
either $\dim f(B_1)\leq 1$ or $\dim f(B_1\cap B_2)=0$. \qed
\end{enumerate}
\end{cor}

In the log canonical setting we also used the
$i=1$ case of the following.
In contrast with (\ref{relkod.klt.nonbir.thm}) and
(\ref{relkod.lc.cor}), its proof is quite difficult and subtle.

\begin{thm}\label{amb-fuj.thm}\cite[3.2]{amb}, \cite[2.39]{fujinobook},
\cite[6.3]{fuj-fund}
 Let $(Y, \sum_iD_i)$  be an embeddable snc pair
and  $f:Y\to X$  a projective morphism.
Let $M$ be a line bundle on $Y$ and assume that
$$
M\sim_{\q,f}K_Y+(\mbox{$f$-semi-ample})+\tsum_i a_i D_i
\qtq{where $0\leq a_i\leq 1$.}
$$
 Then
every associated prime of $R^if_*M$ is  the $f$-image
of an irreducible component of some intersection
$B_1\cap\cdots\cap B_r$ for
some divisors $B_j\subset \sing Y\cup\{D_i: a_i=1\}$. \qed
\end{thm}

\begin{rem}[The analytic case]\label{anal.ref.rem}
 The complex analytic version of
(\ref{relkod.klt.thm}) is proved in  \cite{tak, nak}
but the complex analytic version of
(\ref{amb-fuj.thm}) is not known. 

If $X$ is a complex analytic space,
we can choose $f:Y\to X$ to be projective.
By pushing forward the sequence
$$
0\to M\bigl(-\tsum D_i\bigr)\to M\to M|_{\sum D_i}\to 0
$$
shows that $R^if_*M\cong R^if_*\bigl(M|_{\sum D_i}\bigr)$.
Thus we need to prove the analog of (\ref{amb-fuj.thm})
for the projective morphism $\sum D_i\to X$.
This is in  fact how the proofs of \cite{amb, fujinobook}
work. However, their proofs rely on a global
compactification of $\sum D_i\to X$.
\end{rem}

\begin{say}[A natural double cover]\label{double.cover.trick}
 Let $X^0$ be a  scheme whose singularities are
double nnormal crossing points only. Let $\pi^0:\bar X^0\to X^0$ denote
its normalization with conductors  $D^0\subset X^0$, 
 $\bar D^0 \subset \bar X^0$
and Galois involution $\tau:\bar D^0 \to \bar D^0$.

Take two copies $\bar X^0_1\amalg \bar X^0_2$  and
on $\bar D^0_1\amalg \bar D^0_2$
consider the involution
$$
\sigma(p,q)=\bigl(\tau(q), \tau(p)\bigr).
$$
Note that
$\bigl(\bar D^0_1\textstyle{\amalg} \bar D^0_2\bigr)/\sigma\cong \bar D^0$
but the isomorphism is non-canonical.

We obtain $\tilde X^0$ either as the universal push-out 
of
$$
\bigl(\bar D^0_1\amalg \bar D^0_2\bigr)/\sigma\leftarrow
\bigl(\bar D^0_1\amalg \bar D^0_2\bigr)\into
\bigl(\bar X^0_1\amalg \bar X^0_2\bigr)
$$
or as the spectrum of the preimage of the $\sigma$-invariant part 
$$
\bigl(\o_{\bar D^0_1}+\o_{\bar D^0_2}\bigr)^{\sigma}
\subset \o_{\bar D^0_1}+\o_{\bar D^0_2}
\qtq{in}
\o_{\bar X^0_1}+\o_{\bar X^0_2}. 
$$
Then $\pi^0:\tilde X^0\to X^0$ is an \'etale double cover
and the irreducible components of $\tilde X^0$ are smooth.
The normalization of $\tilde X^0$ is a disjoint union of
two copies of the normalization of $X^0$.

Let  now $X$ be  slc and $j:X^0\into X$ an open subset
with double nc points only and such that 
$X\setminus X^0$ has codimenson $\geq 2$.
Let $\pi^0:\tilde X^0\to X^0$ be as above.
Then $j_*\pi^0_*\o_{ \tilde X^0}$ is  a coherent sheaf of
algebras on $X$. Set
$$
\tilde X:=\spec_X j_*\pi^0_*\o_{\tilde X^0}
$$
with projection $\pi:\tilde X\to X$.

By construction, $\tilde X$ is $S_2$,  $\pi$
 is \'etale in codimension 1 and 
the normalization of $\tilde X$ is a disjoint union of
two copies of the normalization of $X$.
Furthermore, the
irreducible components of $\tilde X$ are smooth
in codimension 1.
\end{say}

\subsection*{Conditions for (\ref{d.two.ss.say.3}.3)}{\ }

These  are reduced to the previous vanishing
theorems using the following duality.

\begin{prop} \label{duality.prop}
Let $f:Y\to X$ be a proper morphism, $Y$ CM, 
$M$ a vector bundle on $Y$  and $x\in X$ a closed point.
Set $n=\dim Y$ and 
let $W\subset Y$ be a subscheme such that
$\supp W=\supp f^{-1}(x)$.
Then there is a natural bilinear pairing
$$
H^i_W\bigl(Y, \omega_Y\otimes M^{-1}\bigr)\times\bigl( R^{n-i}f_* M\bigr)_{x}
\to k(x)
$$
which has no left or right kernel,
where the subscript denotes the stalk at $x$.

In particular, 
if either $H^i_W\bigl(Y, \omega_Y\otimes M^{-1}\bigr)$ or
$\bigl( R^{n-i}f_* M\bigr)_{x}$ 
 is a finite dimensional $k(x)$-vector space
then so is the other and they are dual to each other.
\end{prop}

Proof. Let $mW\subset Y$ be the subscheme defined by
the ideal sheaf $\o_Y(-W)^m$. 
By \cite[II.6]{sga2}
$$
H^i_W\bigl(Y, \omega_Y\otimes M^{-1}\bigr)=
\varinjlim \ext^i_Y\bigl(\o_{mW},\omega_Y\otimes  M^{-1}\bigr).
\eqno{(\ref{duality.prop}.1)}
$$
On the other hand, by the theorem on formal functions,
$$
\bigl( R^{n-i}f_* M\bigr)^{\textstyle \wedge}=
\varprojlim H^{n-i}\bigl(mW,  M|_{mW}\bigr)
\eqno{(\ref{duality.prop}.2)}
$$
where $\wedge$ denotes completion at $x\in X$.

We show below that for every $m$, the groups on the right hand sides of
(\ref{duality.prop}.1--2) are dual to each other.
This gives the  required bilinear pairing
which has no left or right kernel.
\medskip

(Duality using a compactification.)
Let $\bar Y\supset Y$ be a CM compactification
such that $M$ extends to a line bundle $\bar M$ on $\bar Y$.
Since $W$ is disjoint from  $\bar Y\setminus Y$,
$$
\ext^i_Y\bigl(\o_{mW},\omega_Y\otimes  M^{-1}\bigr)=
\ext^i_{\bar Y}\bigl(\o_{mW}, \omega_{\bar Y}\otimes  \bar M^{-1}\bigr)=
\ext^i_{\bar Y}\bigl(\o_{mW}\otimes  \bar M, 
\omega_{\bar Y}\bigr)
$$ and, 
by Serre duality, 
the latter is  dual to 
$$
H^{n-i}\bigl(\bar Y, \o_{mW}\otimes \bar M\bigr)
=
H^{n-i}\bigl(mW, \bar M|_{mW}\bigr)=
H^{n-i}\bigl(mW,   M|_{mW}\bigr).
$$
(It is not known that such $\bar Y$ exists,
but we could have used any compactification and  Grothendieck duality.
For the complex analytic case see \cite{RRV}.)
 
\medskip

(Duality without   compactification.)
This proof works if 
$W$ is an effective Cartier divisor.
In most applications, $X$ is given and $Y$ is a suitable resolution,
hence this assumption is easy to achieve.

We use the local-to-global spectral sequence for $\ext$
$$
H^i\bigl(Y, \sext^j_Y(N,N')\bigr)\Rightarrow \ext^{i+j}_Y(N,N').
$$
Since $\o_{mW}$ has projective dimension 1 as an $\o_Y$-sheaf,
$ \sext^i_Y\bigl(\o_{mW}, \omega_Y\bigr)=0$ for $i\neq 1$
and $ \sext^1_Y\bigl(\o_{mW}, \omega_Y\bigr)=\omega_{mW}$.
Thus for any locally free $\o_{mW}$-sheaf $N$, 
the local-to-global spectral sequence for $\ext^*_Y(N,\omega_Y)$
degenerates and
$$
\begin{array}{rcl}
\ext^i_Y(N,\omega_Y)& = &
H^{i-1}\bigl(Y,  \sext^1_Y(N, \omega_Y)\bigr)\\
&=&H^{i-1}\bigl(Y,  N^{-1}\otimes \omega_{mW}\bigr)\ =\
H^{i-1}\bigl(mW,  N^{-1}\otimes \omega_{mW}\bigr).
\end{array}
$$
Setting $N=\o_{mW}\otimes M$ gives
the  isomorphisms
$$
 \ext^i_Y\bigl(\o_{mW}, \omega_Y\otimes M^{-1}\bigr)=
 H^{i-1}\bigl(mW, M^{-1}\otimes \omega_{mW}\bigr).
$$
Since $mW$ is a projective CM scheme over a field, 
Serre duality gives that the latter group is dual to
$$
H^{n-1-(i-1)}\bigl(mW,  M|_{mW}\bigr)=
H^{n-i}\bigl(mW,   M|_{mW}\bigr). \qed
$$

\medskip

Combining this with (\ref{relkod.klt.thm}) we obtain the following.

\begin{cor} \label{duality.cor}
Let $f:Y\to X$ be a proper morphism, $Y$  smooth.
Let $x\in X$ be a closed point and
assume that  $W:=\supp f^{-1}(x)$ is a Cartier divisor.
Let $L$ be a line bundle on $Y$  such that
$L\sim_{\q,f}(\mbox{$f$-nef})+\Delta$ for some simple normal crossing
 divisor $\Delta$ such  that $\rdown{\Delta}=0$. Then 
  $H^i_W\bigl(Y, L^{-1}\bigr)=0$ for $i<\dim X$. \qed
\end{cor}

There are obvious dual versions of (\ref{relkod.klt.nonbir.thm})
and of (\ref{relkod.lc.cor}).

 \begin{ack} I thank O.~Fujino and S.~Mori for
 useful comments and corrections.
Partial financial support  was provided by  the NSF under grant number 
DMS-0758275.
\end{ack}

\bibliography{refs}

\vskip1cm

\noindent Princeton University, Princeton NJ 08544-1000

\begin{verbatim}kollar@math.princeton.edu\end{verbatim}

\end{document}